\input psfig
\font\titlefont= cmcsc10 at 12pt

\def\R{{\bf R}}
\def\C{{\bf C}}
\def\Z{{\bf Z}}
\def\Q{{\bf Q}}
\font\eighteenbf=cmbx10 scaled\magstep3
\def\vandaag{\number\day\space\ifcase\month\or
 januari\or februari\or  maart\or  april\or mei\or juni\or  juli\or
 augustus\or  september\or  oktober\or november\or  december\or\fi,
\number\year}
\def\h{{h^0}}
\def\k{{k^0}}
\def\ZZ{{\bf Z}}
\def\FF{{\bf F}}
\def\CC{{\bf C}}
\def\RR{{\bf R}}
\def\QQ{{\bf Q}}

\def\zx{{Z^{}_X}}
\def\zf{{Z^{}_F}}
\def\abs#1{|\!|#1|\!|}

\def\sqd{{\sqrt{|\Delta|}}}
\magnification\magstep1
\vglue 1in\centerline{\eighteenbf  Effectivity of Arakelov Divisors and}
\bigskip
\centerline{\eighteenbf  the Theta Divisor of a
Number Field}
\bigskip
\vskip 2pc
\centerline{\titlefont Gerard van der Geer \& Ren\'e Schoof}
\bigskip\bigskip
\centerline{\bf 1. Introduction}
\smallskip
\noindent
In the  well known analogy between the theory of function fields of curves over
finite fields  and the arithmetic of algebraic number fields, the number theoretical
analogue of a divisor on a curve is an Arakelov divisor. In this paper we
introduce the notion of an {\it effective} Arakelov divisor; more precisely, we
attach to every Arakelov divisor $D$ its {\it effectivity},  a real number
between~0 and~1. This notion naturally leads to another quantity associated
to~$D$. This is a positive real number $\h(D)$ which is the arithmetic
analogue of the dimension of the vector space $H^0(D)$ of sections of the
line bundle associated to a divisor~$D$ on an algebraic curve. It can be
interpreted as the logarithm of a value of a theta function. Both notions
can be extended to higher  rank Arakelov bundles. 
\par
In this paper we show that the {\it effectivity} and the numbers $h^0(D)$ behave
in several respects like there traditional geometric analogues. 
We prove an analogue of the fact that $h^0(D)=0$ for any
divisor $D$ on an algebraic curve with ${\rm deg}(D)<0$. We provide evidence for a
conjecture that would imply an analogue of the fact that $h^0(D)\le {\rm
deg}(D)+1$ for divisors $D$ on an algebraic curve with ${\rm deg}(D)\ge 0$.
The Poisson summation formula implies a Riemann-Roch Theorem involving the
numbers~$\h(D)$ and $h^0(\kappa -D)$ with
$\kappa$ the canonical class; it is a special case of Tate's
Riemann-Roch formula. Unfortunately, we do not have a definition of
$h^1(D):=h^0(\kappa -D)$ for an Arakelov divisor without recourse to duality. 
Following K.~Iwasawa~[Iw] or
J.~Tate~[T] one derives in a natural way the finiteness of the class group and the
unit theorem of Dirichlet from this analogue of the Riemann-Roch Theorem, avoiding
the usual arguments from geometry of numbers.

The notion of effectivity naturally  leads to a 
definition of the zeta function of a number field which is closely analoguous to
the zeta function of a curve over a finite field. In this way the
Dedekind zeta function, multiplied by the usual gamma factors,  is recovered as an
integral over the Arakelov class group. 

There is a close connection between Arakelov divisors and certain lattices. The
numbers $h^0(D)$ are closely related to the Hermite constants of these lattices.
The value of $h^0$ on the canonical class is an
invariant of the number field that can be viewed as an analogue of the genus of a
curve.  The quantity $\h(D)$ defines a real analytic function on the Arakelov
divisor class group. Its restriction to the group of Arakelov divisors of degree
${1 \over 2} \deg(\kappa)$  can be viewed as the analogue of the set of rational
points of the theta divisor of an algebraic curve over a finite field. 
\par
It is natural to try to obtain arithmetic analogues of various basic
geometric facts like Clifford's Theorem. As is explained at the end of
section~5, this comes down to studying the behaviour of the function
$h^0$ on a space parametrizing  bundles of rank~2.
\par
We  suggest in the same spirit a definition for an invariant $h^0(L)$ for
a metrized line bundle $L$ on an arithmetic surface. The definition for
$\h$ provided here hints at a further theory and we hope this paper will
stimulate the readers to develop it.
\par
In section~2 we recall some well known facts concerning Arakelov divisors.
In section~3 we introduce the notion of effectivity and the
definition of~$\h$. We apply these to the zeta function of a number field
in section~4. In section~5 we give some estimates on $h^0(D)$. We introduce the
analogue of the genus for number fields in section~6. In section~7 we briefly
discuss a two variable zeta function. Finally, in section~8 and~9, we make some
remarks about the higher dimensional theory.
\smallskip
\noindent
{\bf Acknowledgements.} The first author would like to thank Dick Gross,
Hendrik Lenstra, A.N.\ Parshin and Don Zagier for helpful discussions on
the theme of this paper. We also thank Wieb Bosma for help with Figures 4~and~5 and
Richard Groenwegen for catching some errors in the calculations in an earlier
version of section~6.
\bigskip

\bigskip
\centerline{\bf 2. Arakelov Divisors}
\smallskip
\noindent
The similarity between the class group $Cl(O_F)$ of the ring of integers
of a number field~$F$ and the Jacobian of an irreducible smooth curve is a
particular aspect of the deep analogy between number fields and function fields
of curves. However, the class group classifies isomorphism classes of
line bundles on the {\it affine} scheme ${\rm Spec}(O_F)$ and algebraic geometry
tells us that the Jacobian of a {\it projective} or complete curve is much better
behaved than that of an affine one.  Following S.~Arakelov [A1,A2] one arrives at an
improved analogy via a sort of compactification of the scheme ${\rm Spec}(O_F)$
by taking the archimedean primes of $F$ into account. One thus obtains a
generalization of the class group which is a compact group; it is an extension
of the usual (finite) ideal class group by a real torus. In this section we recall
the definitions and briefly discuss a variant of this theory.
\bigskip
Let $F$ be a number field. An {\it Arakelov divisor} is a formal finite sum
$\sum_Px_P P+\sum_{\sigma}x_{\sigma}\sigma$, where
$P$ runs over the prime ideals of the ring of integers $O_F$ and $\sigma$ runs
over the infinite, or archimedean, primes of the number field~$F$. The
coefficients
$x_P$ are in~$\ZZ$ but the $x_{\sigma}$ are in~$\RR$. The Arakelov
divisors form an additive  group ${\rm Div}(F)$ isomorphic to
$\sum_P\ZZ\times\sum_{\sigma}\RR$. The first sum is infinite but the second
is a real vector space of dimension $r_1+r_2$. Here $r_1$ and $r_2$ denote the
number of real and complex infinite primes, respectively. We have that
$r_1+2r_2=n$ where $n=[F:\QQ]$. The degree ${\rm deg}(D)$ of an Arakelov divisor
$D$ is given by
$$
{\rm deg}(D)=\sum_P{\rm log}(N(P))x_P+\sum_\sigma x_{\sigma}.
$$
The norm of $D$ is given by $N(D)=e^{{\rm deg}(D)}$.
\par
An Arakelov divisor $D=\sum_Px_P P+\sum_\sigma x_{\sigma}\sigma$ is
determined by the associated fractional ideal $I=\prod P^{-x_P}$ and by the
$r_1+r_2$ coefficients
$x_{\sigma}\in\RR$ at the infinite primes. 
For every $f\in F^*$ the principal Arakelov divisor $(f)$ is 
defined by $(f)=\sum_Px_P P+\sum_{\sigma}x_{\sigma}\sigma$, where
$x_P={\rm ord}_P(f)$ and 
$x_{\sigma}=-{\rm log}|\sigma(f)|$ or $-2{\rm log}|\sigma(f)|$ depending on
whether $\sigma$ is real or complex. The ideal associated to $D$ is the
principal fractional ideal $f^{-1}O_F$.
By the product formula we have ${\rm deg}(f)=0$.
The principal Arakelov divisors form a subgroup of the group ${\rm Div}(F)$.
The quotient group is called the Arakelov divisor class group or
Arakelov-Picard group. It is denoted by ${\rm Pic}(F)$. There is an exact
sequence
$$
0\longrightarrow\mu_F\longrightarrow F^*\longrightarrow {\rm
Div}(F)\longrightarrow {\rm Pic}(F)
\longrightarrow 0.
$$
Here $\mu_F$ denotes the group of roots of unity in~$F^*$. 
Since the degree of a principal divisor is zero, the degree map ${\rm
Div}(F)\longrightarrow \RR$ factors through
${\rm Pic}(F)$. For
$d\in\RR$, we denote by  ${\rm Pic}^{(d)}(F)$  the set of divisor classes of
degree~$d$.  Forgetting the infinite components, we obtain a surjective
homomorphism from ${\rm Pic}^{(0)}(F)$ to the ideal class group $Cl(O_F)$ of the
ring of integers ~$O_F$ which fits into an exact sequence 
$$
0\longrightarrow V/\phi(O_F^*)\longrightarrow {\rm
Pic}^{(0)}(F)\longrightarrow Cl(O_F)\longrightarrow 0.
$$
Here $V=\{(x_{\sigma})\in\prod_{\sigma}\RR:\sum_{\sigma} x_{\sigma}=0\}$ and 
$\phi:O_F^*\longrightarrow\prod_{\sigma}\RR$ is the natural map
$O_F^*\longrightarrow {\rm Div}(F)$ followed by the projection on the
infinite components. The vector space $V$ has dimension $r_1+r_2-1$ and
$\phi(O_F^*)$ is a discrete subgroup of~$V$. In section~4 we show 
that the group ${\rm Pic}^{(0)}(F)$ is compact. This statement is
equivalent to Dirichlet's Unit Theorem and the fact that the ideal class group
$Cl(O_F)$ is finite.  The volume of ${\rm Pic}^{(0)}(F)$ is equal to $hR$, where
$h=\#Cl(O_F)$ and
$R$ denotes the regulator of~$F$. 
\par
It is natural to associate a lattice to an Arakelov divisor
$\sum_Px_P P+\sum_{\sigma}x_{\sigma}\sigma$.
When $\sigma$ is real, the
coefficient $x_{\sigma}$ determines a scalar product on $\RR$ by setting
$\abs{1}_{\sigma}^2=e^{-2x_{\sigma}}$. When  $\sigma$ is complex, $x_{\sigma}$
determines a hermitian product on $\CC$ by setting
$\abs{1}_{\sigma}^2=2e^{-x_{\sigma}}$. Taken together, these metrics induce a
metric on the product $\RR^{r_1}\times\CC^{r_2}$ by 
$$
\abs{(z_{\sigma})}_D^2 =\sum_{\sigma}|z_{\sigma}|^2\abs{1}_{\sigma}^2 \, .
$$
We view, as usual, the number field $F$ as a subset of $\RR^{r_1}\times\CC^{r_2}$
via the embeddings $\sigma$.   With these metrics the covolume of
the lattice $I$ is equal to $\sqd/N(D)$, where $\Delta$ denotes the discriminant
of~$F$. 
\par
It is not difficult to see that the classes of two Arakelov divisors $D$~and~$D'$ 
in ${\rm Pic}(F)$ are the same if and only if there is an $O_F$-linear
isomorphism $I\longrightarrow I'$ that is compatible with the metrics on the
associated lattices. Therefore, the group
${\rm Pic}(F)$ parametrizes isometry classes of
lattices with compatible $O_F$-structures. The cosets ${\rm Pic}^{(d)}(F)$
parametrize such lattices of covolume $\sqd e^{-d}$.
\smallskip
For an Arakelov divisor $D$ one defines the
Euler-Poincar\'e characteristic (cf.\ [Sz]):
$$ 
\chi (D) = -  \log ({\hbox{covol}}(I))={\rm deg}(D)-{1\over 2}{\rm log}|\Delta|, 
$$
where the covolume is that of the ideal associated to~$D$, viewed as
a lattice $ I \subset \R^{r_1}\times\CC^{ r_2}$ equipped with the metrics
induced by~$D$. Since $\chi (O_F)= -{1 \over 2} \log |\Delta |$, we have for an
Arakelov divisor $D$ that 
$$
\chi(D)={\rm deg}(D)+\chi(O_F).
$$
\bigskip
\noindent
{\bf Variant.} We may also consider the following variant of Arakelov theory.
Let
$$
W= \prod_{\sigma} \C / \Z (1),
$$
with  $\Z(1)=2\pi i \Z$ and where the product runs over {\it all}
embeddings. Define a  complex conjugation $c$ op $W$ via:
$$
c: z_{\sigma} \mapsto \bar{z}_{\bar{\sigma}}.
$$
We have
$$
\eqalign{ W^c & = \{ (z_{\sigma})_{\sigma} : \bar{z}_{\sigma} =
z_{\bar{\sigma}} \} \cr
&= (\R \times \pi i \Z / 2\pi i \Z)^{r_1} \times (\C / 2\pi i\Z)^{r_2}.
\cr}
$$
Consider Arakelov divisors similar to the ones above but now of the form
$$
D= \sum_P x_P\,  P + \sum_{\sigma} \lambda_{\sigma} \, \sigma,
\qquad \hbox{\rm with } \bar{\lambda}_{\sigma} =
\lambda_{\bar{\sigma}},
$$
where the sum is again over all embeddings  $\sigma$.
The principal divisor $(f)$ of an element  $f\in F^*$ is defined by
$$
(f)= \sum_P {\rm ord}_P(f) \, P + \sum_{\sigma} \log (\sigma(f)).
$$
The resulting class group is denoted by  $\widetilde{\rm Pic}(F)$. We have an exact
 sequence
$$
1 \longrightarrow O_F^* \longrightarrow W^c \longrightarrow 
\widetilde{\rm Pic}(F) \longrightarrow Cl(O_F) \longrightarrow 1.
$$
The real  dimension of
$\widetilde{\rm Pic}(F)$ is $n= [F:\Q]$. Comparison with the usual
Arakelov-Picard group  can be done by mapping the infinite coefficients
$\lambda_{\sigma}$ to their real parts. This induces a surjective  map
$\widetilde{\rm Pic}(F) \to {\rm Pic}(F)$ whose kernel is isomorphic to the group
$\prod_{\sigma} (i\RR/2\pi i\ZZ)^c
\cong  (\Z/2\Z)^{r_1} \times (\RR/2\pi\ZZ)^{r_2}$ modulo the group $\mu_F$ of the
roots of unity of~$F$.
\par
\bigskip
\noindent
{\bf Remark.} Instead of the full ring of integers one can also consider
orders $A \subset O_K$ and consider these as the analogues of singular
curves. Define
$$
\partial_{A/\Z}^{-1} := {\rm Hom}_{\Z}(A,\Z)=\{ x \in K : {\rm Tr}(xA)
\subset \Z \}.
$$
If $A$ is Gorenstein this is a locally free $A$-module of rank~1 and we  can define
the canonical divisor to be the ideal $(\partial_{A/ \Z}^{-1})^{-1}$ together with
the standard metrics. For instance,  if  $A=\Z [\alpha]= \Z[x]/(f(x))$ then the
canonical divisor is $(f^{\prime}(\alpha))$  by Euler's identity. We
can now develop the theory for these orders instead for $O_F$.  But we
can also change the metrics on $\partial_{F/\Q}$ which means a change
of model at the infinite places. The interpretation of this is
less clear.

\bigskip
\centerline{\bf 3. Effectivity and an Analogue of the Theta Divisor}
\smallskip
\noindent
A divisor $D=\sum_Px_PP$ of a smooth complete absolutely irreducible algebraic
curve~$X$ over a field~$k$ is called {\it effective}  when $n_P\ge 0$ for all
points~$P$. The vector space $H^0(D)$ of sections of the associated line bundle is
defined as
$$
H^0(D)=\{f\in k(X)^*: \hbox{$(f)+D$ is effective}\}\cup\{0\}.
$$
Here $(f)$ denotes the principal divisor associated to a non-zero function 
$f\in k(X)$. 
The usual generalization~[Sz] of the space $H^0(D)$ to Arakelov divisors 
$D=\sum_Px_PP+\sum_{\sigma} x_{\sigma}\sigma$ of a number field~$F$ is given by
$$
\eqalign{
H^0(D)&=\{f\in F^*: \hbox{all coefficients of the divisor $(f)+D$ are
non-negative}\}\cup\{0\},\cr
&= \{ f \in I \colon \|f\|_{\sigma} \leq 1\,
{\hbox { for all infinite primes $\sigma$} }  \}. \cr}
$$
Here $I=\prod_PP^{-x_P}$ denotes the ideal associated to~$D$.
The set $H^0(D)$ is the intersection of a lattice and a compact set.
Therefore it is finite and one then puts
$$
\h(D):= {\rm log}(\# H^0(D)).
$$
This is not very satisfactory. We  introduce here a new
notion of $h^0(D)$ for Arakelov divisors $D$.
First we introduce the notion of {\it effectivity} of an Arakelov divisor.
Let $D$ be an Arakelov divisor of a number field~$F$. We view $F$ as a subset
of the Euclidean space $\RR^{r_1}\times\CC^{r_2}$ via its real and complex
embeddings~$\sigma$. For an element $f\in F$ we write $\abs{f}_D$ for
$\abs{(\sigma(f))}_D$.
\par
We define the {\it effectivity} $e(D)$ of an Arakelov divisor $D$ by
$$
e(D)=\cases{0,&if $O_F\not\subset I$;\cr \exp({-\pi\abs{1}_{D}^2}),& if
$O_F\subset I$.\cr} 
$$
We have that $0\leq e(D)<1$. Explicitly, for $D=\sum_Px_P P
+\sum_{\sigma}x_{\sigma}\sigma$ we have that $e(D)=0$ whenever $x_P<0$ for some
prime~$P$. If $x_P\ge 0$ for all~$P$, we have that $$
e(D)=\exp({-\pi\abs{1}_D^2})=\exp({-\pi\!\!\!\!\sum_{\sigma \,\,\,\rm real}\!\!\!\!
e^{-2x_{\sigma}}-\pi\!\!\!\!\!\!\!\sum_{\sigma \,\,\,\rm complex}\!\!\!\!\!\!
2e^{-x_{\sigma}}}).
$$
The effectivity of $D$ is close to~1 when each
$x_\sigma$ is large. If one of the $x_{\sigma}$ becomes negative however, the
effectivity of
$D$ tends doubly exponentially fast to~$0$. For instance, for $F=\QQ$ the
effectivity of the Arakelov divisor $D_x$ with finite part $\ZZ$ and infinite
coordinate $x_{\sigma}=x$ is given by the function $e(D_x)=e^{-\pi e^{-2x}}$, see
Fig.\ 1.
\smallskip
\centerline{\psfig{figure=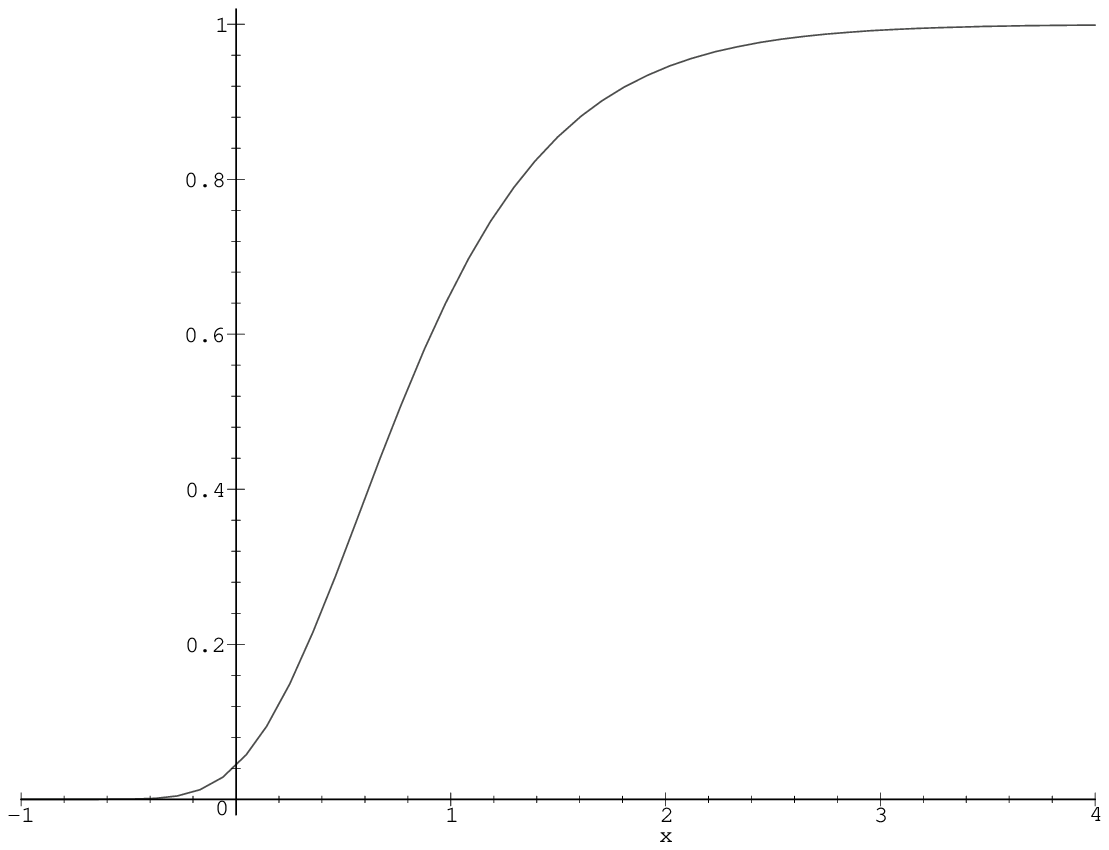,width=6cm}}
\par\bigskip\par
\smallskip
\centerline{$\scriptstyle{\rm \,\,  Fig.1.\,
The\,\, function\,\,}
e(D){\,\,\rm  for\,\,}  \bf Q .$}
\medskip\noindent
We can use the notion of effectivity to define the analogue of the dimension
$\h(D)$ of the vector space $H^0(D)$ of sections of the line bundle associated to
a divisor $D$ of a curve. It seems natural to put for an Arakelov divisor~$D$
$$
H^0(D)=\{f\in F^*:e((f)+D)>0\}\cup\{0\}.
$$
The effectivity of $(f)+D$ is  positive if and only if $f\in I$.
Here $I$ denotes the fractional ideal associated to~$D$. Therefore  $H^0(D)$ is
equal to the infinite group $I$. The function  $e$ attaches a weight 
(`the effectivity') to  non-zero functions $f \in H^0(D)$ viewed as sections of 
the bundle $O(D)$ via  $f \mapsto e((f)+D)$. We consider $H^0(D)$ together with the
effectivity as the analogue of the geometric $H^0(D)$.
 To measure its size, we
weight the divisors $(f)+D$ to which the elements $f$ give rise with their
effectivity 
$$
e((f)+D)=e^{-\pi\abs{1}_{(f)+D}^2}= 
e^{-\pi\abs{f}_{D}^2}.
$$
When we count elements $f$ rather than the ideals they generate and add
the element $0\in I$ to the sum, we obtain the analogue of the order of the
group~$H^0(D)$:
$$
\k(D)= \sum_{f\in I} e^{-\pi\abs{f}_{D}^2}.
$$
See also~[Mo, 1.4]. The analogue of the dimension of $H^0(D)$ is given by
$$
\h(D)={\rm log}\left(\sum_{f\in I} e^{-\pi\abs{f}_{D}^2}\right),
$$
which we call the {\it size} of $H^0(D)$. Since two Arakelov divisors in the
same class in ${\rm Pic}(F)$ have isometric associated lattices, the function
$\h(D)$ only depends on the class $[D]$ of~$D$ in ${\rm Pic}(F)$ and we may write
$\h[D]$. This is a function on the Picard group.

\medskip
There is an analogue of the  Riemann-Roch Theorem for the numbers $\h(D)$.
We define the {\it canonical divisor} $\kappa$ as the Arakelov divisor whose
ideal part is the inverse of the {\it
different} $\partial$ of $F$ and whose infinite components are all zero. We have
that
$N(\partial)=|\Delta|$, so that ${\rm deg}(\kappa)=\log|\Delta|$.
Therefore
${1\over 2}{\rm log}|\Delta|$ may be seen as the analogue of the quantity
${g-1}$ that occurs in the Riemann-Roch formula for curves of genus~$g$.
\medskip

\proclaim Proposition 1. (Riemann-Roch) Let $F$ be a number field with discriminant
$\Delta$ and let $D$ be an Arakelov divisor. Then
$$
h^0(D) - h^0(\kappa-D) = {\rm deg}(D) -{1\over 2}{\rm log}|\Delta|.
$$
\par
\noindent {\bf Proof.} This is Hecke's functional equation for the theta
function. The lattices associated to $D$ and $\kappa-D$ are $\ZZ$-dual to one
another and the formula follows from an application of the Poisson summation
formula:
$$
\sum_{f\in I} e^{-\pi\abs{f}_{D}^2} = {{N(D)}\over{\sqd}}
\sum_{f\in \partial I^{-1}} e^{-\pi\abs{f}_{\kappa-D}^2}.
$$
\bigskip
This Riemann-Roch theorem is a special case of Tate's~[T]. Just like Tate had many
choices for his zeta functions, we had a choice for our effectivity function. 
We chose the function $e^{-\pi \abs{x}^2}$, because it gives rise to a symmetric
form of the Riemann-Roch theorem and because it leads to the functional equation
of the Dedekind zeta function. This is analogous to the  geometric case, where
there is a unique function $h^0$. 
  We do not have a definition of $h^1(D)$ without recourse to duality.
However, recently A.~Borisov~[Bo] proposed a direct definition of $h^1(D)$.  

\medskip
Recall that in the geometric case the Jacobian of a curve~$X$ of genus~$g$ together
with its theta divisor $\Theta\subset {\rm Pic}^{(g-1)}(X)$ determines the
curve and recall that the geometry of the curve (e.g.\ existence of linear
systems) can be read off from the theta divisor. In particular, Riemann showed
that the singularities of $\Theta$ determine the linear systems on~$X$ of
degree~$g-1$: 
$$
\h(D)={\rm ord}_{[D]}(\Theta),
$$
where ${\rm ord}_{[D]}(\Theta)$ is the multiplicity of $\Theta$
at the point $[D] \in {\rm Pic}^{g-1}(X)$.

Let $d={1\over 2}{\rm
log}|\Delta|$. We view the restriction of the function $\h$ to ${\rm Pic}^{(d)}(F)$
as the analogue of the theta divisor~$\Theta$.   The function $\h$ is a
real analytic function on the space ${\rm Pic}^{(d)}(F)$. It should be possible
to reconstruct the arithmetic of the number field $F$ from ${\rm Pic}^{(d)}$
together with this function.
\smallskip
We give some numerical examples.

\smallskip\noindent
{\bf Example 1.} For $F=\QQ$ the function $h^0$ looks as follows. Since $\ZZ$ has
unique factorization, the degree map ${\rm Pic}(\QQ)\longrightarrow \RR$ is an
isomorphism. To $x\in\RR$ corresponds the divisor $D_x$ that has associated ideal
equal to $\ZZ$ and infinite coordinate $x_{\sigma}=x$. We have that $\h(D_x)= {\rm
log}(\sum_{n\in\ZZ} e^{-\pi n^2 e^{-2x}})$.  This function tends very rapidly to
zero when $x$ becomes negative. For instance, for $x=-3$ its value is smaller than
$10^{-500}$. The Riemann-Roch Theorem says in this case that $h^0(D_x)
-h^0(D_{-x})=x$ for all $x\in\RR$. 
\smallskip
\centerline{\psfig{figure=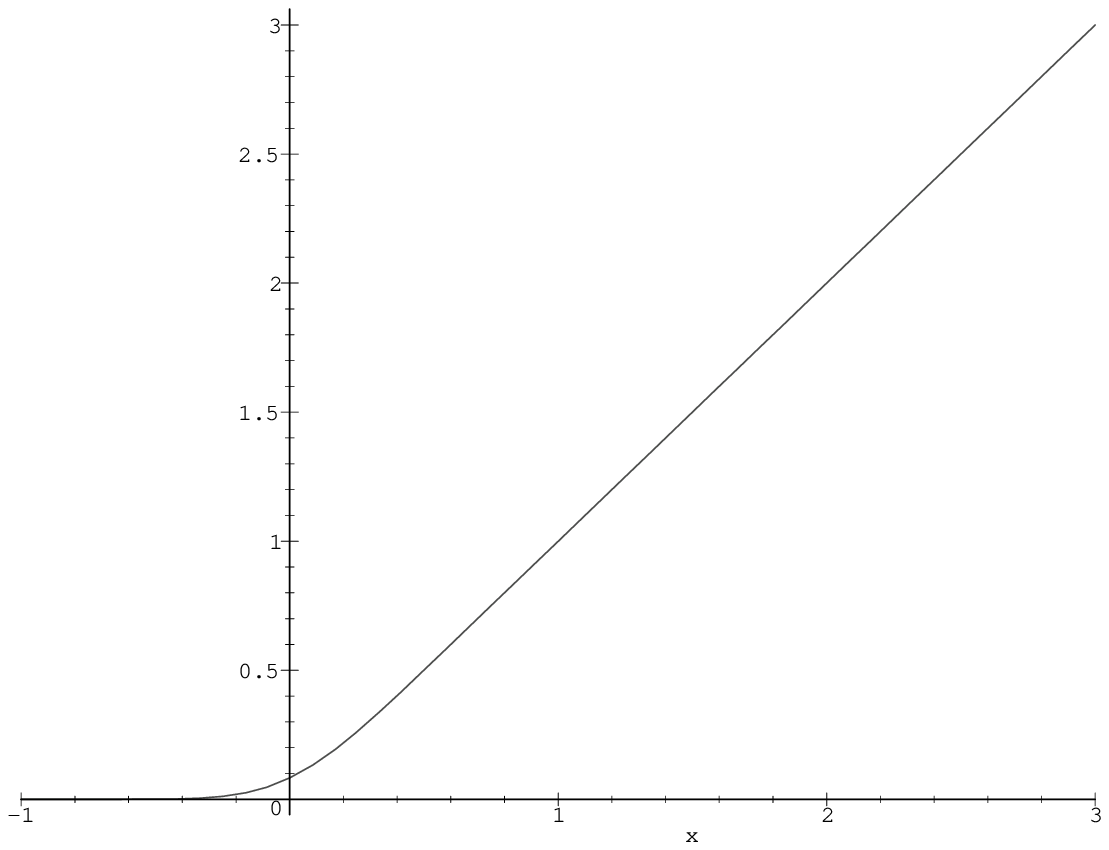,width=5cm}}
\par\bigskip
\par
\smallskip
\centerline{$\scriptstyle{\rm \,\,  Fig.2.\,
The\,\, function\,\,}
\h(D){\,\,\rm  for\,\,}  \bf Q .$}
\medskip\noindent
Since $D_x$ is the Arakelov divisor $x\, \sigma$ where $\sigma$ is the unique
infinite prime of~$\QQ$, the function $\h(D_x)$ is the analogue of
the function $\h(nP)={\rm dim}\,H^0(nP)$, where $P$ is a point on the projective
line ${\bf P}^1$. In that case one has that $\h(nP)={\rm
max}(0,n+1)$. 
\medskip
\noindent
{\bf Example 2.} For a real quadratic field of class number~1, the Picard group
is an extension of $\RR$ by ${\rm Pic}^{(0)}(F)$. The group ${\rm
Pic}^{(0)}(F)$ is isomorphic to $\RR/R\ZZ$ where $R$ is the regulator
of~$F$.  We take $F=\QQ(\sqrt{41})$ and we plot the function
$\h(D_x)$ on ${\rm Pic}^{(0)}(F)$. Here $D_x$ denotes the divisor whose ideal part
is the ring of integers
$O_F$ of $F$ and whose infinite part has coordinates~$x$ and~$-x$. It
is periodic modulo $R=4.15912713462618001310854497\ldots$ Note the relatively big
maximum for $x=0$. This phenomenon is the analogue of the geometric fact that for a
divisor $D$ of degree $0$ on a curve we have $h^0(D)=0$ unless $D$ is the trivial
divisor, for which $h^0(D)=1$. 
\smallskip
\centerline{\psfig{figure=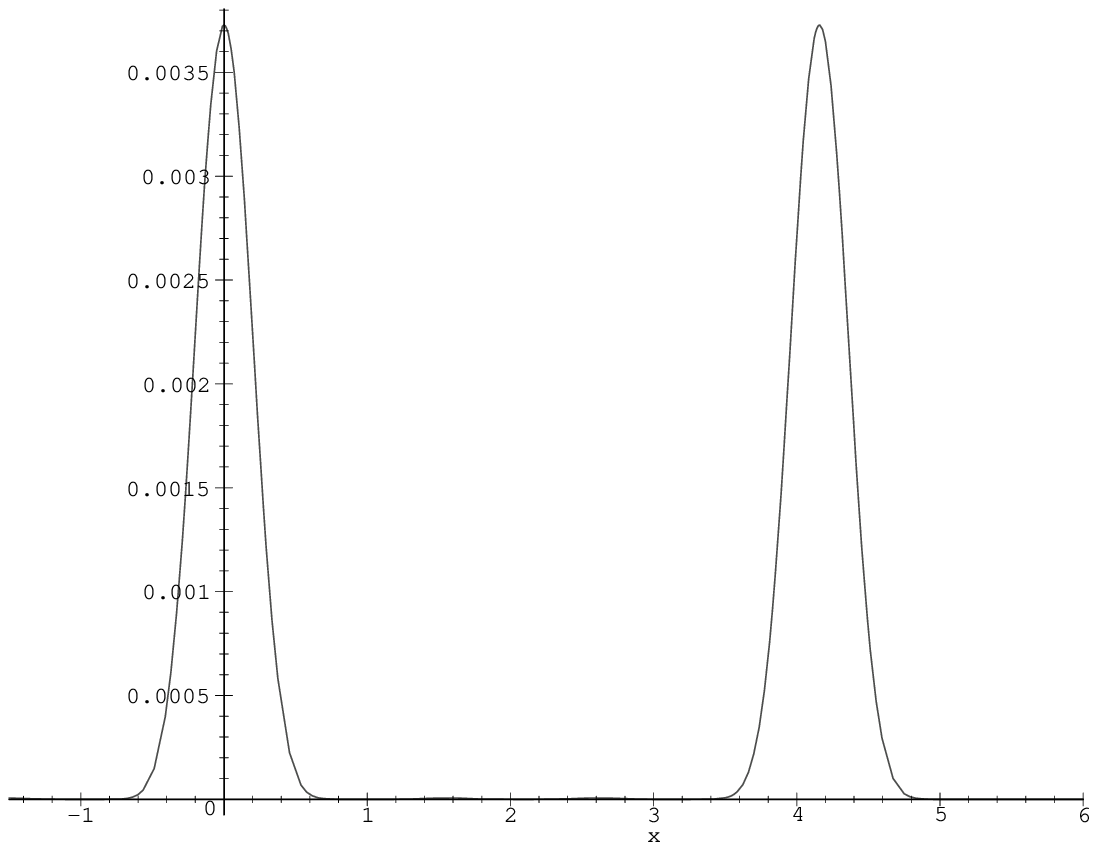,width=8cm}}
\par\smallskip\par
\bigskip
\centerline{$\scriptstyle{\rm \,\,  Fig.3.\,
The\,\, function\,\,}
\h(D){\,\,\rm  on\,\,} {\rm Pic}^{(0)}(\QQ(\sqrt{41}))$.}
\smallskip\noindent
\noindent
{\bf Example 3.} We now take $F=\QQ(\sqrt{73})$. The class number of $F$ is~1
and the group ${\rm Pic}(F)$ is isomorphic to a cylinder. For every $d\in\RR$, the
coset ${\rm Pic}^{(d)}$ of classes of degree~$d$ is a circle whose circumference 
is equal to the regulator $R=7.666690419258287747402075701\ldots$ of~$F$. The
classes of ${\rm Pic}^{(d)}$ are represented by the divisors $D_x$ whose
integral part is $O_F$ and whose infinite components are $d/2-x$ and $d/2+x$
respectively. We depict the functions
$\h(D_x)$ restricted to ${\rm Pic}^{(d)}(F)$  for $d={i\over{10}}{\rm
log}|\Delta|$ with $i=0,1,\ldots,9$. The Riemann-Roch theorem says that the
functions for $i$ and $10-i$  are translates of one another by
${{|5-i|}\over{10}}{\rm log}(73)$.
\bigskip
\centerline{\psfig{figure=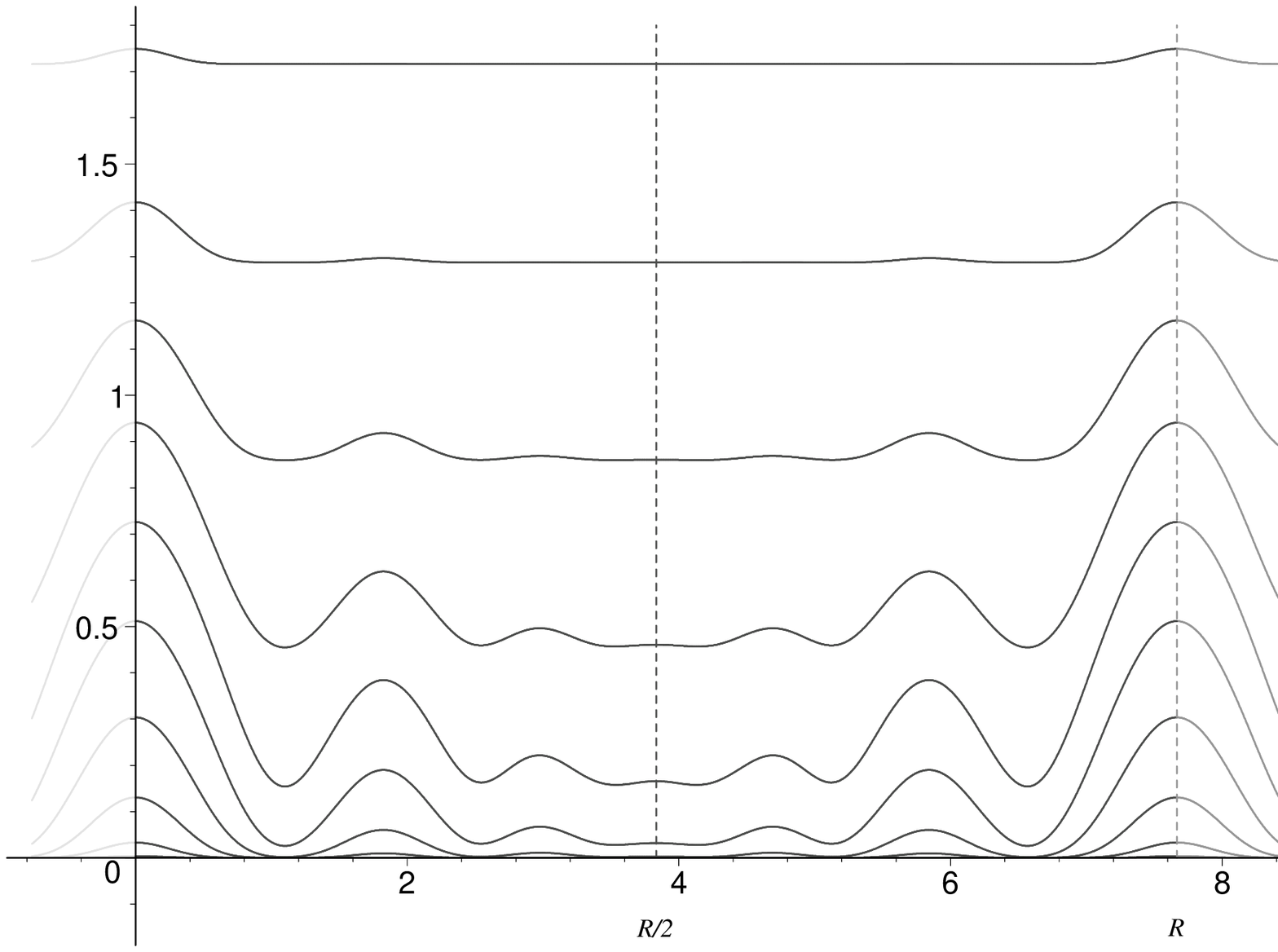,width=7cm}}
\par
\bigskip
\par
\bigskip
\centerline{$\scriptstyle{\rm \,\,  Fig.4.\,
The\,\, function\,\,}
h^0(x){\,\,\rm  for\,\,}  {\bf Q}(\sqrt{73}) .$}

\bigskip\bigskip
\centerline{\bf 4. Zeta functions}
\smallskip
\noindent
In this section we recover the zeta function $\zeta_F(s)$ of a number
field $F$ as a certain integral of the ``effectivity" function on the group of
Arakelov divisors. It is a natural adaptation of the definition of the zeta
function of a curve over a finite field. Following 
Iwasawa~[Iw] and Tate~[Ta], we use the Riemann-Roch Theorem to prove that the
topological group ${\rm Pic}^{(0)}(F)$ is compact. This gives a proof of the
finiteness of the class group and of Dirichlet's Unit Theorem that only makes use
of the functional equation of the Theta function and is not based on the usual
techniques from geometry of numbers. The geometric analogue of this result is the
theorem that for a curve over a finite field the group ${\rm Pic}^{(0)}$ is finite.
The usual proof of this fact exploits the Riemann-Roch theorem and the fact that
the number of {\it effective} divisors of fixed degree is finite. From our point
of view, the proof by Iwasawa and Tate is a natural generalization of this argument.
\smallskip

We briefly discuss the zeta function of an algebraic curve over a finite field.
Let $X$ be an absolutely irreducible complete smooth algebraic curve of
genus~$g$ over $\FF_q$. We denote the group of $\FF_q$-rational divisors by ${\rm
Div}(X)$. The degree of a point is the degree of 
its residue field over~$\FF_q$.
The degree of a divisor $D=\sum_Pn_P\cdot P\in{\rm Div}(X)$ is given by
${\rm deg}(D)=\sum_Pn_P{\rm deg}(P)$. We let $N(D)=q^{{\rm deg}(D)}$ denote the
norm of~$D$.   The zeta function $\zx(s)$ of $X$ is defined by
$$
\zx(s)=\sum_{D\ge 0}N(D)^{-s},\qquad\qquad(s\in\CC,\,\,{\rm Re}(s)>1).
$$
Here the sum  runs over the {\it effective} divisors $D$ of $X$. It converges
absolutely when ${\rm Re}\,(s)>1$. In order to analyze the zeta function, one
considers the Picard group
${\rm Pic}(X)$ of~$X$. This is the group 
${\rm Div}(X)$ modulo the group of $\FF_q$-rational  principal divisors 
$\{(f):f\in\FF_q(X)^*\}$. Since the degree of a principal divisor is zero, the
degree ${\rm deg}(D)$ and the norm $N(D)$ only depend on the class $[D]$ of
$D$ in the Picard group. By
${\rm Pic}^{(d)}(X)$ we denote the divisor classes in the Picard group that have
degree~$d$.

We rewrite the sum, by first summing over divisor classes $[D]\in{\rm Pic}(X)$
and then counting the effective divisors in each~$[D]$, i.e.\ counting $|D|$, the
projectivized vector space of sections
$$
H^0(D) = \{f\in \FF_q(X)^*:{\rm div}(f)+D\ge 0\}\cup\{0\}
$$
of dimension~$h^0(D)$. We let $\k(D)=\#H^0(D)$. Since the norm $N(D)$ only depends
on the divisor class $[D]$ of~$D$, we write $N[D]$. We have    
$$
\zx(s)=\sum_{[D]\in {\rm Pic}(X)}\#\{D'\in[D]:\hbox{$D'$
effective}\}N[D]^{-s} =\sum_{[D]\in {\rm
Pic}(X)}{{\k(D)-1}\over{q-1}}N[D]^{-s}. 
$$
The Riemann-Roch Theorem implies that $\k(D)=q^{{\rm deg}(D)-g+1}$ when ${\rm
deg}(D)>2g-2$ and it relates $\k(D)$ to $\k(\kappa-D)$. Here $\kappa$ denotes the
canonical divisor of~$X$. These facts imply that one can sum the series
explicitly and that $\zx(s)$ admits a meromorphic continuation to all of~$\C$.
Moreover, $\zx(s)$ has a simple pole at $s=1$ with residue given by
$$
\mathop{\rm Res}\limits_{s=1}\zx(s)={{\#{\rm Pic}^{(0)}(X)}\over
{(q-1)q^{g-1}{\rm log}(q)}}.
$$
In this way one can actually {\it prove} that ${\rm Pic}^{(0)}(X)$ is finite.
In some sense, this proof is an analytic version
of the usual argument that there are only finitely many
effective divisors of fixed degree. The fact that the zeta function converges
for $s>1$ is a stronger statement that generalizes better to number fields.
\smallskip

Next we turn to the arithmetic situation.
First we give a definition of the zeta function associated to a number
field~$F$ that is natural from the point of view of section~3. Rather than summing
$N(I)^{-s}$ over all non-zero ideals $I$ of $O_F$, we take ---just like we did for
curves--- the sum over all effective Arakelov divisors~$D$. More precisely,
$$
\zf(s)=\int_{{\rm Div}(F)}N(D)^{-s}d\mu_{\rm eff},
$$
where $\mu_{\rm eff}$ denotes the measure on ${\rm Div}(F)$ that weights the
Arakelov divisors with their effectivity. To see that this integral converges
absolutely for $s\in\CC$ with ${\rm Re}(s)>1$, we split the integral into a
product of an infinite sum and a multiple integral. Writing $J$ for the {\it
inverse} of the ideal associated to~$D$ and $t_{\sigma}$ for $e^{-x_{\sigma}}$,
where the
$x_{\sigma}$ denote the infinite components of~$D$, we find
$$
\eqalign{\zf(s)&=\int_{{\rm Div}(F)}e(D)N(D)^{-s}dD,\cr
&=\sum_{0\not=J\subset O_F}N(J)^{-s}\int_{t_{\sigma}}
\left(\prod_{\sigma}t_{\sigma}^{s}\right) \exp({-\pi\!\!\sum_{\sigma\,\,\rm
real}t_{\sigma}^2-\pi\!\!\!\!\!\!\sum_{\sigma\,\,\rm
complex}\!\!\!\!2t_{\sigma}})\prod_{\sigma}{{dt_{\sigma}}\over{t_{\sigma}}},\cr
&=\left(2\pi^{-s/2}\Gamma(s/2)\right)^{r_1}\left(
(2\pi)^{-s}\Gamma(s)\right)^{r_2}\sum_{0\not=J\subset O_F}N(J)^{-s}.\cr}
$$
We see that our zeta function is precisely the Dedekind
zeta function multiplied by the usual
gamma factors. Therefore the integral converges absolutely for $s\in\CC$
with ${\rm Re}(s)>1$.  
\smallskip
\noindent
This way of writing $Z_F(s)$ may serve to motivate the definition
of the effectivity~$e$. As we remarked above, a priori other definitions of the
effectivity are possible, and this leads to modified zeta functions similar to the
integrals in Tate's thesis [T]. Our choice is the one made by Iwasawa in his
1952 letter to Dieudonn\'e~[Iw]. In this letter Iwasawa also showed that one
can establish the compactness of ${\rm Pic}^{(0)}(F)$ by adapting the computations
with the zeta function of a curve over a  finite field that we indicated above. We
briefly sketch his arguments.
\par
Just as we did for curves, we write the zeta function as a repeated
integral
$$
\zf(s)=\int_{{\rm
Pic}(F)}N[D]^{-s}\left(\int_{[D]}e^{-\pi\abs{1}_{D}^2}dD\right)d[D].
$$
The divisors in the coset~$[D]$ have the form $(f)+D$. Only the ones for which
$f$ is contained in $I$, where $I$ is the ideal associated to $D$, have
non-zero effectivity. Therefore
$$\eqalign{
\int_{[D]}e^{-\pi\abs{1}_{D}^2}dD&=\sum_{(0)\neq (f)\subset
I}e^{-\pi\abs{1}_{(f)+D}^2} =\sum_{(0)\neq (f)\subset
I}e^{-\pi\abs{f}_{D}^2},\cr
&={1\over w}\left(-1+\sum_{f\in
I}e^{-\pi\abs{f}_{D}^2}\right),\cr}
$$
where $w$ is the number of roots of $1$ in $O_F$. This gives the following
expression for
$Z_F(s)$:
$$
\zf(s)={1\over w}\int_{{\rm
Pic}(F)}(k^0[D]-1)N[D]^{-s}d[D].
$$
An application of Proposition 1 (i.e. the arithmetic Riemann-Roch formula)
easily implies the folowing
$$\eqalign{
w\zf(s)=\int_{{[D]\in {\rm Pic}(X)}\atop{N[D]<\sqd}}
(k^0(D)-1)N[D]^{-s}&d[D]+
\int_{{[D]\in {\rm Pic}(X)}\atop{N[D]\le\sqd}}
(k^0(D)-1){{N[D]^{s-1}}\over{\sqrt{|\Delta |}}^{(2s-1)}}d[D]\cr
&\qquad \qquad \qquad \qquad+{{{\rm vol}({\rm
Pic}^{(0)}(F))}\over{s(s-1)\sqd^{s}}}.\cr}
$$
For $s\in\RR$, $s>1$, all three summands are  positive. 
Since $Z_F(s)$ converges for $s>1$, substituting any such $s$ gives therefore an
upper bound for the volume of ${\rm Pic}^{(0)}(F)$. As explained in section~2 
the fact that ${\rm Pic}^{(0)}(F)$ has finite volume, implies Dirichlet's Unit
Theorem and the finiteness of the ideal class group $Cl_F$. 
\par
In contrast to the situation for curves over finite fields, this time the
effective Arakelov divisors $D$ of negative degree or, equivalently, with
$N(D)<1$, contribute to the integrals. However their contribution is very small
because they have been weighted with their effectivity.  We estimate the integrals 
in the next section. This  estimate is also used there to deduce the meromorphic
continuation of $Z_F(s)$.
\par
\bigskip\bigskip
\centerline{{\bf 5. Estimates for} $h^0(D)$}
\bigskip
\noindent
In this section we give an estimate on $h^0(D)$ and discuss an analogue of the
inequality $h^0(D) \leq \deg(D) +1$ for divisors $D$ with ${\rm deg}(D)\ge 0$.
We also describe the relation between $h^0(D)$ and the Hermite constant of the
lattice associated to~$D$. 
\smallskip
\noindent
\proclaim Proposition 2.   Let $F$ be a number field of degree~$n$. Let $D$ be an 
Arakelov divisor $D$  of~$F$ with ${\rm deg}(D)\le{1\over 2}{\rm log}|\Delta|$ and
let $f_0\in I$ be the shortest non-zero vector in the lattice $I$ associated
to~$D$. Then 
$$
k^0(D)-1 \leq \beta e^{-\pi \abs{f_0}_{D}^2}.
$$
for some constant $\beta$ depending only on the field~$F$.
\par
\smallskip
\noindent
{\bf Proof.} Let $D$ be a divisor with $N(D)\leq \sqrt{|\Delta|}$, or equivalently,
$\deg(D)\leq {1\over 2}{\rm log}|\Delta|$. We define a positive real number $u$ by
$u=({1\over 2}\log(|\Delta|)-\deg(D))/n$ with $n=[F:\Q]$.
Define a new divisor $D^{\prime}$ of degree ${1\over 2}{\rm log}|\Delta|$ by
$$
D^{\prime}= D + \sum_{\sigma \, {\rm real}} u\, \sigma  + \sum_{\sigma \, {\rm
complex}} 2u \, \sigma.
$$
For $0\neq f\in I$, the ideal associated to $D$, we have $\abs{f}^2_D -
\abs{f}^2_{D^{\prime}} =
\abs{f}_{D}^2 (1-e^{-2u})$, hence
$$
\eqalign{e^{-\pi \abs{f}_{D}^2}&=e^{-\pi \abs{f}_{D^{\prime}}^2}\cdot 
e^{-\pi(\abs{f}_{D}^2 - \abs{f}^2_{D^{\prime}})} \cr
&\leq e^{-\pi \abs{f}_{D^{\prime}}^2}\cdot 
e^{-\pi(\abs{f_0}_{D}^2 - \abs{f_0}^2_{D^{\prime}})}, }
$$
where $f_0$ is a shortest non-zero vector for $D$, hence also for $D^{\prime}$.
We get
$$
k^0(D)-1 \leq (k^0(D^{\prime})-1)\cdot e^{-\pi \abs{f_0}_{D}^2}\cdot
e^{\pi \abs{f_0}_{D^{\prime}}^2}.
$$
In section~4 we showed that the cosets ${\rm Pic}^{(d)}(F)$
are compact. The functions 
$k^0(D^{\prime})-1$ and $e^{\pi \abs{f_0}_{D^{\prime}}^2}$ are continuous on the
coset of divisor classes of degree $d={{1\over 2}{\rm log}|\Delta|}$, hence are
bounded. This implies the Proposition.
\medskip
\proclaim Corollary 1.  Let $F$ be a number field of degree~$n$. Let $D$ be an 
Arakelov divisor $D$  of~$F$ with ${\rm deg}(D)\le{1\over 2}{\rm log}|\Delta|$.
Then 
$$
h^0(D)<k^0(D)-1 \leq \beta e^{-\pi ne^{-{2\over n}{\rm deg}(D)}}
$$
for some constant $\beta$ depending only on the field~$F$.

\smallskip\noindent{\bf Proof.} Let $x_{\sigma}$ denote the infinite components
of~$D$ and put $t_{\sigma}=e^{x_{\sigma}}$.  Let $0\not=f\in I$ be a non-zero
vector in~$I$. By the geometric-arithmetic mean inequality we have
$$\eqalign{
\abs{f}^2&=\sum_{\sigma\,\,\rm real}t^{-2}_{\sigma}|\sigma(f)|^2+
\sum_{\sigma\,\,\rm complex}2t_{\sigma}^{-1}|\sigma(f)|^2\cr
&\ge n(|N(f)|^2\prod_{\sigma}t_{\sigma}^{-2})^{1/n}\ge 
n(|N(I)|^{2}\prod_{\sigma}t_{\sigma}^{-2})^{1/n}=nN(D)^{-2/n}.\cr}
$$
It follows that  $e^{-\pi \abs{f}_{D}^2}$ is bounded from above by $e^{-\pi\, n
e^{-{2\over n}\deg(D)}}$. Proposition~2 now implies the result.
\medskip

Corollary~1 is the  analogue of the geometric fact that $H^0(D)=0$
whenever $D$ is a divisor on a curve with ${\rm deg}(D)<0$.
Indeed, if  ${\rm deg}(D)$ becomes negative, then the proposition implies
that $\h(D)$ tends doubly exponentially fast to zero. 

The estimate of Corollary~1 leads to the meromorphic continuation of the zeta
function~$Z_F(s)$.  Indeed, since there exists a constant $\beta$  only
depending on the number field~$F$ so that $0\le k^0(D)-1\le
\beta e^{-\pi n N(D)^{-2/n}}$ whenever $N(D)<\sqd$,
the two integrals in the last expression for $wZ_f(s)$ in section~4 converge rapidly
to functions that are holomorphic in~$s\in \CC$. This implies that the zeta function
extends to a meromorphic function on~$\CC$. The function $|\Delta|^{s/2}Z_F(s)$ is
invariant under the substitution $s \mapsto 1-s$. It is not difficult to see that
the residue of
$Z_F(s)$ at the pole in $s=1$ is given by
$$
{{{\rm vol}({\rm Pic}^{(0)})}\over{w\sqd}}.
$$
\proclaim Corollary 2. Let $F$ be a number field of degree~$n$ and let $w=\#\mu_F$
denote the number of roots of unity in~$F$. Then there is a constant $\beta$
depending only on~$F$ so that
$$
\left|{{{\rm log}\,({1\over w}h^0(D))}\over{{\rm
covol}(I)^{2/n}}}+\pi\gamma(I)\right|<\beta e^{2{\rm deg}(D)/n},\qquad
\hbox{for all divisors $D$ with ${\rm deg}(D)<0$.}
$$
Here $I$ denotes the ideal associated to~$D$ and $\gamma(I)$ denotes the
{\it Hermite constant} of the lattice associated to~$I$. In other words $\gamma(I)$
is the square of the length of the shortest non-zero vector in~$I$ divided by ${\rm
covol}(I)^{2/n}$.

\smallskip\noindent{\bf Proof.} Let $f_0$ denote the shortest non-zero vector
of~$I$. We obviously have that $h^0(D)\ge {\rm log}(1+ we^{-\pi\abs{f_0}^2})$.
Combining this with the inequality of
Prop.2 and using the compactness of ${\rm Pic}^{(0)}(F)$ to bound $\abs{f_0}^2$
from below, we find that
$$
\beta' e^{-\pi\abs{f_0}^2}\le h^0(D)\le \beta'' e^{-\pi\abs{f_0}^2}
$$
for certain $\beta',\beta''>0$. Dividing these quantities by~$w$, taking the
logarithm and dividing by~${\rm covol}(I)^{2/n}$ easily implies the corollary.

\medskip
The Hermite constant only depends on the lattice modulo homothety.
The corollary says that the function $-{1\over{\pi}}{\rm log}({1\over w}h^0(D)){\rm
covol}(I)^{-2/n}$ approaches the Hermite constant $\gamma(I)$ of $I$ when ${\rm
deg}(D)$ tends to $-\infty$. In practice the convergence is rather fast. We give a
numerical example.

\smallskip\noindent{\bf Example.} As in section~3, we consider $F=\QQ(\sqrt{73})$
and the function $h^0$ on ${\rm Pic}(F)$. We depict a graph of the function
$$
B^0(d,x)=-{{{\rm log}({1\over 2}\h(D_x))}\over{2\pi{\rm exp}(-d)}}
$$
for $d=0$ and for $d=-{1\over 2}{\rm log}|\Delta|$. Here $D_x$
denotes the divisor whose integral part is $O_F$ and whose infinite coordinate
are ${d\over2}+x$ and ${d\over2}-x$ respectively. Its divisor class is contained in
${\rm Pic}^{(d)}(F)$. The lattice associated to $D_x$ is denoted by~$I_x$. 
Its covolume is equal  to $\sqrt{|\Delta|}{\rm exp}(-d)$. It
follows from Corollary~2 that $B^0(d,x)$ tends to
$\gamma(I_x){{\sqrt{73}}\over 2}$ as $d$ tends to~$-\infty$.  We see that the
graphs for $d=0$ and $d=-{1\over 2}{\rm log}|\Delta|$ are extremely close.  The
graph for $d=-{1\over 2}{\rm log}|\Delta|$ visibly outdoes the one for
$d=0$ only near its local maxima. \par
 
\centerline{\psfig{figure=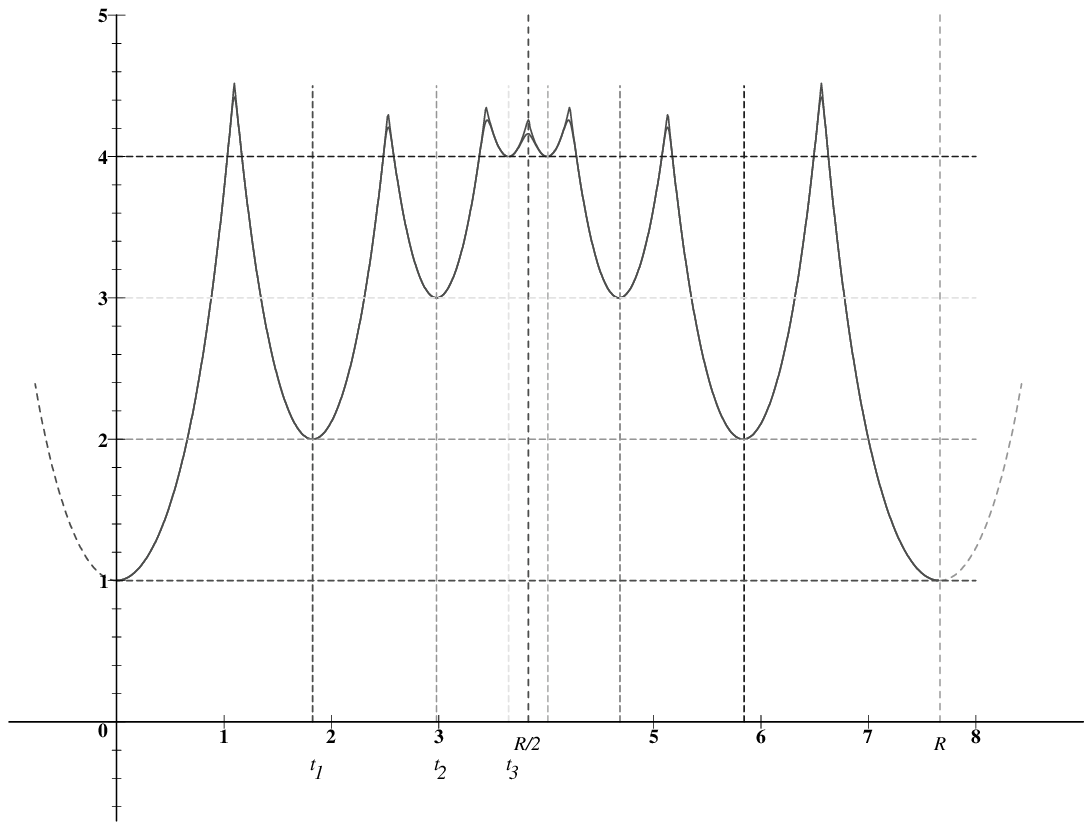,width=12cm}}
\par\bigskip\par\smallskip
\centerline{$\scriptstyle{\rm \,\,  Fig.5.\,
The\,\, function\,\,}
B^0(d,x){\,\,\rm  for\,\,}  {\bf Q}(\sqrt{73}) .$}
\medskip\noindent
For every $d\in\RR$, the function $B^0(d,x)$ is periodic modulo
the regulator~$R={\rm log}(1068+125\sqrt{73})=7.666690419258287747402075701\ldots$
and  symmetric with respect to~$0$. The graph reflects properties of the lattice
$O_F\subset\RR\times\RR$. The maxima of the function $h^0(D_x)$, and hence the 
minima of $B^0(x)$, are attained at ${1\over 2}{\rm log}|f/\overline{f}|$ where $f$
is very near one of the successive minima $1$, $(9+\sqrt{73})/2$, $17+2\sqrt{73}$,
$(77+9\sqrt{73})/2$ or one of their conjugates. Here $\overline{f}$ denotes the ${\rm
Gal}(F/\QQ)$-conjugate of~$f$. Because of our normalization, the value of
$B^0(d,x)$ in these points is approximately equal to $|f\overline{f}|$.
\bigskip
In the geometric case the function $h^0$ on ${\rm Pic}^0(X)$ assumes the value~$1$
in the trivial class and is $0$ elsewhere. We conjecture a similar behaviour in the
arithmetic case. Roughly speaking, we expect that the function $h^0$ assumes a
pronounced maximal value in or very close to the trivial class $[O_F]$.
Computations suggest that this does indeed usually happen. However, if there is
a unit in $O_F^*$ of infinite order, all of whose absolute values are relatively
close to~1,  the function $h^0$ may assume its maximum value rather far away
from~$[O_F]$. In this case however it seems that $h^0[O_F]$ is nevertheless
rather close to the maximum value of~$h^0$.

If the number field $F$ admits many automorphisms, we dare be more precise.
\noindent
\proclaim  Conjecture.  Let $F$ be a number field that is Galois over $\bf Q$ or
over an imaginary quadratic number field. Then the function
$h^0$ on
${\rm Pic}^0$  assumes its maximum in the trivial class $O_F$.

\smallskip\noindent
It is easy to see that under these conditions the function $h^0$ has a local maximum
in the trivial class  $[O_F]
\in {\rm Pic}(F)$. The conjecture has been proved by P.~Francini~[Fr] for quadratic
number fields.

The analogue of the geometric fact $\dim H^0(D) \le
\deg(D) + 1$ if  $\deg(D) \geq 0$ follows from this conjecture: 
\smallskip

\proclaim
Proposition 3. If the function $h^0(D)$ on ${\rm Pic}^0(F)$  assumes its maximum in
the origin $[O_F]$, then for every $D$ with $\deg(D) \geq 0$ we have $h^0(D) \leq
\deg(D) + h^0(O_F)$.
\par
\smallskip
\noindent
{\bf Proof.} We set $u= \deg(D)/n$ and we define a divisor $D^{\prime}$ of degree $0$
by 
$$
D^{\prime} = D- \sum_{\sigma \, \,
{\rm real}}u \, \sigma - \sum_{\sigma \, \, {\rm complex}}2u \, \sigma .
$$
Then $\deg(D^{\prime})= 0$. By a term-by-term comparison for $D$ and $D^{\prime}$
we get $h^0(\kappa-D)\le  h^0(\kappa -D^{\prime})$.
But then we have
$$
\eqalign{
h^0(D)-\deg (D) &= h^0(\kappa -D) - \log(\sqrt{|\Delta|}) \qquad {\hbox{\rm (by
Riemann-Roch)}}\cr 
&\leq h^0(\kappa - D^{\prime}) -  \log(\sqrt{|\Delta|})\cr 
&=h^0(D^{\prime})\qquad {\hbox {\rm (by Riemann-Roch)}}\cr
&\leq h^0(O_F) \qquad {\hbox{\rm (by assumption).}}
}
$$
This proves the Proposition.
\medskip
Finally we mention Clifford's theorem. This classical result is a statement about
the function $h^0(D)$ for divisors $D$ of an algebraic curve of genus~$g$. It says
that for every divisor $D$ with $0\le {\rm deg}(D)\le 2g-2$ one has that $h^0(D)\le
{1\over 2}{\rm deg}(D)+1$. We conjecture that in the arithmetic case the function
$h^0(D)$ behaves in a similar way. In the proof of Clifford's Theorem
one estimates the sum $h^0(D)+h^0(\kappa-D)$. This suggests that in the
arithmetic case one should study the orthogonal sum of the two lattices  associated
to $D$ and $\kappa-D$ and show that the ``size" of this rank 2 lattice is maximal
when $D=O_F$.
\bigskip\bigskip

\vbox{\centerline{\bf 6. A new invariant}
\smallskip
\noindent
By evaluating the function $\k$ at the trivial Arakelov divisor $O_F$
with trivial metrics, we obtain a new invariant of a number field.}
\smallskip
\noindent
{\bf Definition.} The invariant $\eta$  of a number field $F$ is
defined by
$$
\eta(F):= k^0(O_F)= \sum_{x \in \partial^{-1}}e^{-\pi \abs{x}_{\rm triv}^2 }.
$$
Since, by the Riemann-Roch theorem, we have that
$$
\h(\kappa)=\h(O_F) + {1\over2}{\rm log}|\Delta| = \exp(\eta(F)\sqrt{|\Delta|}),
$$
the $\eta$-invariant is directly related to $\h(\kappa)$. The latter should be
viewed as the arithmetic analogue of the {\it genus} of an algebraic curve.

It turns out that the invariant $\eta(F)$ has interesting properties. First we
introduce  a `period':
$$
\omega = {\pi^{1/4} \over \Gamma(3/4)}=1.086434811213308014575316121\ldots
$$
\medskip
\noindent
\proclaim Proposition 4. We have $\eta(\Q)= \omega$. If $F$ is a totally real
number field or a CM-field of degree $n$ then  $\eta(F)= \omega^n \cdot x$
where $x$ is an algebraic number lying in an abelian extension of $\Q(i)$.
\par
\smallskip
\noindent
{\bf Proof.} Let $\theta(\tau)$ be the theta function
$\theta(\tau)=\sum_{n\in\ZZ}e^{\pi in^2\tau}$, ($\tau\in\CC$, ${\rm Im}(\tau)>0$).
The modular function $E_4(\tau)/\theta(\tau)^8$  with $E_4$ the
Eisenstein series of weight $4$, is a rational function on the theta group
$\Gamma_{\vartheta}$, hence assumes by the theory of complex multiplication a
rational value at $\tau=i$, which is $3/4$ as a calculation shows.  Let
$\omega_0$ be the period of the elliptic curve $y^2=4x^3-4x$ defined by
$$
\omega_0=\int_1^{\infty} {dx \over \sqrt{4x^3-4x}}= \int_0^1 {dt \over
\sqrt{1-t^4}}= {1 \over 4}B(1/4,1/2)={\Gamma(1/4)^2\over 4\sqrt{2\pi}},
$$
where $B$ is the Beta function. Since this elliptic curve has multiplication by
$\Z[i]$, the quotient
$\pi^4 E_4(i)/\omega_0^4$ is rational  ($=48$) and the first
statement follows by using the distribution relation satisfied by the
gamma function. For the second statement consider the quotient
$$
{{\sum_{x\in O_F}e^{\pi i \tau\abs{x}^2}}\over{\theta(\tau)^n}}.
$$
This is a modular function with rational coefficients with respect to some
congruence subgroup $\Gamma_0(N)$. Hence by the theory of complex multiplication, 
its value at $\tau=i$ is an algebraic number lying in an abelian extension of
$\Q(i)$. 
\medskip
\noindent
{\bf Examples.} By machine calculation one finds heuristically:
$$
\eqalign{
\eta(\Q(i))&= \omega^2\cdot {2+\sqrt{2}\over 4}, \cr
\eta(\Q(\sqrt{-3}))&= \omega^2 \cdot \left({2+\sqrt{3}\over 4\sqrt{3}}\right)^{1
\over 4},\cr
\eta(\Q(\zeta_7+\zeta_7^{-1}))&= \omega^3 {{
7+3\sqrt{7}+3\sqrt{2\sqrt{7}}}\over{28}}
,\cr
\eta(\Q(\zeta_5))&=\omega^4\cdot {23+\sqrt{5} \over 20} \sqrt{{1+
\sqrt{5}\over 10}}.\cr }
$$
\bigskip\bigskip

\centerline{\bf 7. The Two-variable Zeta Function}
\smallskip
\noindent
In this section we briefly discuss an alternative expression for the
zeta function $Z_X(s)$ of a curve $X$ over~$\FF_q$ which is suggested by the
expression for $(q-1)Z_X(s)$ obtained in section~4. It is related to the two
variable zeta function of Pellikaan~[P]. We also describe the analogue for number
fields.
\noindent
\proclaim Proposition 5. The function
$$
\zeta_X(s,t)=\sum_{[D] \in {\rm Pic}(X)} q^{s\, h^0(D)+t\, h^1(D)}
$$
converges for complex $s,t$ satisfying ${\rm Re}(s)<0,\, {\rm Re}(t) <0$. It can
be continued to a meromorphic function on the whole $(s,t)$-plane. Its restriction
to the line
$s+t=1$ is equal to $(q-1)q^{(g-1)s}Z_X(s)$.
\par\smallskip\noindent
{\bf Proof.} By the Riemann-Roch Theorem, the function
consists of a finite sum $\sum'$ and two infinite parts:
$$
\sum_{[D], \deg(D)>2g-2}q^{s\, (\deg(D)+1-g)} + \sum_{[D],
\deg(D)<0}q^{t\,(-\deg(D)+g-1)}.
$$
Summing the two series, we find
$$
\zeta_X(s,t)=\sum{{}^\prime} + \#{\rm Pic}^{(0)}(X)\left({q^{sg} \over 1-q^s}+
{q^{tg}\over 1-q^t}\right).
$$
This shows that the function admits a meromorphic continuation.
Using the Riemann-Roch Theorem to eliminate the terms $h^1(D)$ in the
finite sum, one finds that the restriction to the line $s+t=1$ is equal to
$(q-1)q^{(g-1)s}Z_X(s)$ as required.
\medskip
Note that the line $s+t=1$ lies entirely outside the domain of convergence.
The following analogue of this result for a number
field~$F$ can be proved in a similar way. \par
\smallskip
\noindent
\proclaim Proposition 6. The function
$$
\zeta_F(s,t)=\int_{{\rm Pic}(F)} e^{s\, h^0(D)+t\, h^1(D)}d[D]
$$
converges for complex $s,t$ satisfying ${\rm Re}(s)<0,\, {\rm Re}(t) <0$. It can
be continued to a meromorphic function on the whole $(s,t)$-plane. Its restriction
to the line
$s+t=1$ is equal to $w {\sqrt{|\Delta|}}^{s/2} \, Z_F(s)$.
\par
For instance, for $F=\QQ$, the function
$$
\zeta_{\bf Q}(s,t)=\int_0^{\infty}\Theta(x)^s\Theta(1/x)^t{{dx}\over x}
$$
converges for ${\rm Re}(s)<0,\, {\rm Re}(t) <0$. Here
$\Theta(x)=\sum_{n\in\ZZ}e^{-\pi n^2x^2}$ is a slight modification of the
usual theta function. The restriction of the analytic continuation to the line
$s+t=1$ is equal to the Riemann zeta function times a gamma factor.
\bigskip\bigskip

\centerline{\bf 8. The higher rank case.}
\smallskip
\noindent
We can extend the definition of Section~3 to higher rank bundles. Let $F$ be a
number field and consider projective $O_F$-modules $M$ of rank $r$ together with 
hermitian metrics at the infinite places. By a theorem of Steinitz a
projective  $O_F$-module of rank~$r$ is isomorphic to  $O_F^{r-1} \oplus I$ for some
ideal $O_F$-ideal $I$.  We define an  Arakelov bundle of rank $r$ to be a projective
$O_F$-module $M$ of rank $r$ together with a hermitian metric on $M\otimes_{\sigma}
\C$ at the infinite places $\sigma$. This defines a metric $\abs{x}_M^2$ on
$M\otimes_{\ZZ}{\RR}$ as it did for fractional ideals.
\par
For such an Arakelov bundle we can define an effectivity on the module  $H^0(M)=M$ of
sections by putting
$$
e(x)=\exp(-\pi\abs{x}_M^2)\qquad {\rm for}\,\,\, x \in M.
$$
and define the size $h^0(M)$  of $H^0(M)$ by
$$
h^0(M)= \log(\sum_{x \in M}  e(x)).
$$
We then have as before a Riemann-Roch theorem. Let 
$$
\chi(M)= -\log {\rm covol}(M).
$$
An easy induction yields $\chi(M)=  \deg(M) + r\chi(O_F)$, where the degree of $M$
is the degree of the Arakelov line bundle  ${\rm det}(M)=\wedge^r M$. Then the
Riemann-Roch theorem is:
$$
h^0(M) -h^0(\kappa\otimes M^{\vee}) = \chi(M),
$$
where $M^{\vee}$ is the dual of $M$ with respect to the Trace map. This is again a
consequence of the Poisson summation formula. 
\par
In order to have reasonable moduli we consider only {\sl admissible}
bundles. This means that the metrics at all infinite   places  $\sigma$ are
obtained by applying an element of  ${\rm SL}(r,O_F)$ to the standard metric
$\sum z_i\bar{z}_i$ on $M\otimes_{\sigma} \CC$.
\par
For simplicity we assume that the projective $O_F$-module $\det(M)$ is actually
free. We  fix a trivialization  $\det (M)\cong O_F$. At the infinite places we put a
(special) hermitian metric on $M\otimes_{\sigma} \C$. Since the admissable metrics
are defined to be the transforms under ${\rm SL}(r,F_{\sigma})$ of the standard
metric, the data at infinity amount to a point of  
$$ {\rm SL}(M) \backslash
\prod_{\sigma} {\rm SL}(M\otimes F_{\sigma})/ {\cal K}, 
$$ 
where ${\cal
K}$ is the stabilizer of a fixed allowable metric at all places $\sigma$. This
quotient can be written as a quotient under ${\rm SL}(M)$ of  
$$  \prod_{\sigma \,\, {\rm real}} {\rm SL}(r,\R)/ {\rm O}(r)  \times
\prod_{\sigma\,\, {\rm complex}} {\rm SL}(r,\C)/ {\rm SU}(r). 
$$
For example, for $r=2$ we find the upper half plane for the real
$\sigma$ and the hyperbolic upper half space  ${\rm SL}(2, \C)/
{\rm SU}(2)$ for the complex infinite places $\sigma$. 
In the particular case that $F$ is totally real and $r=2$ we find the
Hilbert modular varieties associated to $F$. 
\par
As in the geometric case one can now study for a fixed Arakelov
line bundle $L$ of degree $d$ and for varying bundles $M$ of rank $r$ with
$\chi(M)=0$, the function $\Psi_L= h^0(M\otimes L)$ on the moduli space
of Arakelov bundles of fixed rank $r$ with $\chi(M)=0$. This can be viewed as the
analogue of the theta divisors introduced in the geometric case.
\par
\bigskip\bigskip

\centerline{\bf 9. Some Remarks on Higher Dimensions}
\smallskip
\noindent
We finish with some remarks about generalizations. The first remark
concerns the definition of $h^0(L)$ for a metrized line bundle on an arithmetic
surface.
\par
Let $X$ be a smooth projective geometrically irreducible curve over $F$ and assume
that $X$ extends to a semi-stable model ${\cal X}$ over $O_F$. Moreover, we assume
that we are given probability measures $\mu_{\sigma}$ on all $X_{\sigma}(\CC)$. We
consider metrized line bundles $L$ on ${\cal X}$ which are provided with hermitian
metrics at all primes $\sigma$ such that their curvature forms are multiples of
$d\mu_{\sigma}$. One can associate to $L$ the cohomology modules $H^0({\cal X},L)$
and $H^1({\cal X},L)$. These are finitely generated $O_F$-modules. 
It is not possible to define good metrics on them, but Faltings defined a good
metric on the determinant of the cohomology, cf. [F], p.\ 394. We propose to define
an effectivity on $H^0(L)$ as follows. For $s\in H^0(L)$ 
and for each infinite prime $\sigma$ the norm $\abs{s}_{\sigma}$ is defined on
$X_{\sigma}(\CC)$; the divisor of $s$ is of the form $D(s)=
D_f+\sum_{\sigma}x_{\sigma}F_{\sigma}$ with $D_f$ a divisor on ${\cal X}$, with
$F_{\sigma}$ the fibre over $\sigma$ and
$$
x_{\sigma}=-\int_{X_{\sigma}(\CC)}\log\abs{s}_{\sigma}^2d\mu_{\sigma}.
$$
We define the effectivity of $s$ by $e(s)=e(D(s))$ with
$$
e(D(s))={\rm exp}(- \pi \!\!\!\sum_{\sigma\,\, {\rm real } } e^{-2x_{\sigma}}- \pi
\!\!\!\!\!\!\sum_{\sigma \, \, {\rm complex } } 2 e^{-x_{\sigma}})
$$
and the size of $H^0(L)$ by
$$
h^0(L)=\log(\!\!\!\!\sum_{s \in H^0(X,L)}\!\!\!\! e(D(s))).
$$ 
Note that for the trivial line bundle $L$ we get $h^0(O_F)$. Although we do not
have a definition of $h^1(L)$, and we can define $h^2(L)$ only via duality
$h^2(L)=h^0(L^{-1}\otimes \omega_{\cal X})$, one could test whether this definition
is reasonable for suitable very ample line bundles $L$. Then $h^1(L)$ and $h^2(L)$
should be exponentially small, hence our $h^0(L)$ should be close to the Faltings
invariant $\chi(L)$.
\par
As a final remark we point out that a good  notion of effectivity for codimension
$2$ cycles in the sense of Gillet-Soul\'e (cf.\ [G-S]) on an arithmetic surface might
yield a way to write the Hasse-Weil zeta function as an integral over a Chow group.
\bigskip
\vbox{\centerline{\bf References}
\bigskip
\item{[A1]} S.\ Y.\ Arakelov: An intersection theory for divisors on an arithmetic
surface. {\sl Izv.\ Akad.\ Nauk.\ \bf 38} (1974), p.\ 1179--1192. (cf.\ {\sl Math.\
USSR Izvestija \bf 8 } (1974)).
\smallskip
\item{[A2]} S.\ Y.\ Arakelov: Theory of intersections on the arithmetic surface. In
{\sl Proc.\ Int.\ Cong.\ of Math.\ Vancouver}, (1974), p.\ 405--408.
\smallskip
\item{[Bo]} A.\ Borisov: Convolution structures and arithmetic cohomology.
Electronic preprint math.AG/9807151 27 July 1998 at {\tt
http://xxx.lanl.gov} 
\smallskip
\item{[F]} G.~Faltings: Calculus on arithmetic surfaces. {\sl
Ann.\ of Math.\ \bf 119} (1984), p.\ 387--424.
\smallskip
\item{[Fr]} P.\ Francini: The function $h^0$ for quadratic number fields, in
preparation.
\smallskip
\item{[G-S]} H.\ Gillet, C.\ Soul\'e: Intersection sur les vari\'et\'es
d'Arakelov. {\sl C.R.\ Acad.\ Sc.\ \bf 299}, (1984), p.\ 563--566.
\smallskip
\item{[Gr]} R.\ P.\ Groenewegen: The size function for number fields,
Doctoraalscriptie, Universiteit van Amsterdam 1999.
\smallskip
\item{[Iw]} K.\ Iwasawa: Letter to Dieudonn\'e, April 8, 1952, in: N.\ Kurokawa and
T.\ Sunuda: {\sl Zeta Functions in Geometry}, Advanced Studies in Pure Math. {\bf
21} (1992), 445--450.
\smallskip
\item{[Mo]} M.\ Morishita: Integral representations of unramified
Galois groups and matrix divisors over number fields, {\sl Osaka J. Math.} {\bf
32} (1995), 565--576.
\smallskip
\item{[P]} R.\ Pellikaan: On special divisors and the two variable
zeta function of algebraic curves over finite fields. In {\sl Arithmetic,
geometry and coding theory} (Luminy, 1993), 175--184, de Gruyter, Berlin,
1996. 
\smallskip
\item{[Sz1]} L.\ Szpiro: Pr\'esentation de la th\'eorie d'Arakelov. In ``Current
Trends in Arithmetical Algebraic Geometry", {\sl Contemporary Mathematics \bf 67},
AMS, Providence RI (1985), p.\ 279--293. 
\smallskip
\item{[Sz2]} L.\ Szpiro: Degr\'es, intersections, hauteurs. In: S\'eminaire sur
les pinceaux arithm\'e\-ti\-ques: La conjecture de Mordell. {\sl Ast\'erisque \bf
127} (1985), p.\ 11--28.
\smallskip
\item{[T]} J.\ Tate: Fourier Analysis and Hecke's Zeta-functions. Thesis
Princeton 1950. (Printed in: J.W.S.\ Cassels, A.\ Fr\"ohlich: {\sl Algebraic Number
Theory}. Thompson, Washington D.C.\ 1967.)
}
\bigskip
\settabs3 \columns
\+G.\ van der Geer  &&R.\ Schoof\cr
\+Faculteit
Wiskunde en Informatica &&Dipartimento di Matematica\cr
\+Universiteit van
Amsterdam &&II Universit\`a di Roma \cr
\+Plantage Muidergracht 24&&Via Fontanile di Carcaricola \cr
\+1018 TV Amsterdam
&&00133 Roma \cr
\+The Netherlands &&Italy \cr
\+{\tt geer@wins.uva.nl} &&{\tt schoof@wins.uva.nl} \cr
\bye